\documentclass[11pt]{amsart}

\usepackage{todonotes}
\usepackage[T1]{fontenc}
\usepackage[utf8]{inputenc}
\usepackage{lmodern}
\usepackage{microtype}
\usepackage{amsmath,amssymb,amsthm,mathtools}
\usepackage{enumitem}
\usepackage{booktabs}
\usepackage[colorlinks=true,linkcolor=blue,citecolor=blue,urlcolor=blue]{hyperref}

\newtheorem{theorem}{Theorem}[section]
\newtheorem{lemma}[theorem]{Lemma}

\newtheorem{corollary}[theorem]{Corollary}
\theoremstyle{remark}

\newcommand{\Tthree}{T^{3}}
\newcommand{\lk}{\operatorname{lk}}

\newcommand{\defi}{\operatorname{def}}
\newcommand{\normal}[1]{\left\langle\!\left\langle #1\right\rangle\!\right\rangle}
\newcommand{\Q}{\mathbb{Q}}
\newcommand{\Z}{\mathbb{Z}}

\title[No 13-vertex triangulation of the three-torus]
{No 13-vertex triangulation of the three-torus}

\author{Geunho Lim}
\address{Kangnam University, 40 Gangnam-ro, Giheung-gu, Yongin-si, Gyeonggi-do, 16979, Republic of Korea}
\curraddr{}
\email{limg@kangnam.ac.kr}


\begin{document}

\begin{abstract}
K\"uhnel and Lassmann constructed a neighborly triangulation of the three-torus on 15 vertices, and Lutz conjectured that it is vertex-minimal.
The standard lower-bound arguments exclude at most twelve vertices but do not decide the 13- and 14-vertex cases.
We exclude the 13-vertex case.
Our proof uses the Novik--Swartz $g_2$-inequality and a simplicial-complement reduction to a nine-vertex model of the punctured three-torus.
By analyzing fundamental groups of small simplicial complexes and applying Alexander duality in a vertex link, we obtain a contradiction from the deficiency of $\Z^3$.
\end{abstract}

\maketitle

\section{Introduction}
Understanding minimal triangulations of manifolds is a classical problem. The two-torus has a unique vertex minimal triangulation, the Cs\'{a}sz\'{a}r torus on seven vertices.
K\"uhnel and Lassmann constructed a neighborly triangulation of the three-dimensional torus $\Tthree=S^1\times S^1\times S^1$ on 15 vertices \cite{KuhnelLassmann}.  Lutz subsequently conjectured that this triangulation is vertex-minimal
and unique among triangulations with $f$-vector $(15,105,180,90)$
\cite[Conjecture 22]{LutzSurvey}.

The Murai--Novik lower bound for the $g_2$-number of a normal pseudomanifold in terms of the minimum number of generators of its fundamental group gives $f_0(K)\ge 13$ for every triangulation $K$ of $\Tthree$; see \cite{NovikSwartz}.
Thus the first unresolved case is $f_0(K)=13$.  The purpose of this paper is to exclude that case.

\begin{theorem}\label{thm:main}
There is no simplicial triangulation of $\Tthree$ on 13 vertices.
\end{theorem}

Combining Theorem~\ref{thm:main} with the Murai--Novik inequality and the K\"uhnel--Lassmann construction gives the following immediate consequence.

\begin{corollary}\label{cor:bounds}
If $v(\Tthree)$ denotes the minimum number of vertices in a simplicial triangulation of $\Tthree$, then $14\le v(\Tthree)\le 15$.
\end{corollary}

\subsection*{Proof sketch.}
Suppose that $K$ is a 13-vertex triangulation of $\Tthree$.
The Murai--Novik inequality implies that its one-skeleton has at most six missing edges.
Deleting the vertices of a tetrahedral facet $F$ gives a nine-vertex induced subcomplex $C$ with $\pi_1(C)\cong\Z^3$.
The edge bound and the classification~\cite{GovcEtAl} of fundamental groups of complexes on at most eight vertices show that $\lk_C(x)$ is connected and non-simply connected for every vertex $x$ of $C$.
Choosing $x$ to be a vertex in $C$ that is contained in a tetrahedron in $K$ that shares a triangle with $F$, Alexander duality in the two-sphere $\lk_K(x)$ shows that $\pi_1(\lk_C(x))\cong\Z$.
Consequently, $\Z^3$ is obtained from the fundamental group of an eight-vertex complex by adding a single relation.
The classification implies that the fundamental group of this eight-vertex complex has a presentation of deficiency at least two, yielding a presentation of $\Z^3$ of positive deficiency.
This contradicts the exact deficiency $\defi(\Z^3)=0$.

\subsection*{Use of AI}
The authors used ChatGPT 5.6-Sol to assist with typesetting, language editing, proofreading, and locating references.
The authors take full responsibility for all content of this paper.

\section{Preliminaries}

All simplicial complexes are finite.  For a simplicial complex $K$ on vertex set $V(K)$ and a subset $U\subseteq V(K)$, the induced subcomplex on $U$ is denoted by $K[U]:=\{\sigma\in K:\ \sigma\subseteq U\}$.  If $x$ is a vertex, then $K-x:=K[V(K)\setminus\{x\}]$.  The cone on $K$ with apex $x$ is denoted $x*K:=K\cup \{\sigma\cup\{x\}:\ \sigma\in K\}$.
Let $f_i(K)$ denote the number of $i$-dimensional faces in $K$.

In dimension three, every simplicial triangulation of a topological manifold
is a combinatorial triangulation.  In particular, the link of every vertex
in a triangulated closed $3$-manifold is a triangulated $2$-sphere.

For a three-dimensional simplicial complex $K$, set
\[
g_2(K):=f_1(K)-4f_0(K)+10.
\]
If $K$ is a connected normal three-pseudomanifold, Murai and Novik proved
\begin{equation}\label{eq:NS}
 g_2(K)\ge \binom{5}{2}m(\pi_1(K)),
\end{equation}
where $m(G)$ denotes the minimum number of generators of the group $G$ \cite[Theorem 1.1]{NovikSwartz}.

\subsection{The simplicial complement of a face}

The following standard fact is included because the explicit deformation retraction is useful later.

\begin{lemma}\label{lem:complement}
Let $K$ be a simplicial complex, let $\sigma\in K$, put $W=V(K)\setminus\sigma$, and let $C=K[W]$.  Then $|C|$ is a strong deformation retract of $|K|\setminus |\sigma|$.

If, in addition, $K$ triangulates a connected closed three-manifold $M$ and $\sigma$ is a tetrahedral facet, then
\[
\pi_1(C)\cong \pi_1(M).
\]
\end{lemma}

\begin{proof}
Let $p\in |K|\setminus|\sigma|$, and write $\lambda_v(p)$ for its barycentric coordinates.  Put
\[
s(p):=\sum_{w\in W}\lambda_w(p)>0.
\]
For $0\le t\le1$, define barycentric coordinates
\[
\lambda_u(H_t(p))
 =\frac{(1-t)\lambda_u(p)}{1-t(1-s(p))}
 \quad (u\in\sigma),
\]
and
\[
\lambda_w(H_t(p))
 =\frac{\lambda_w(p)}{1-t(1-s(p))}
 \quad (w\in W).
\]
The denominator is positive, the coordinates sum to one, and the support of $H_t(p)$ is contained in the support of $p$.  Hence the homotopy stays in the simplex containing $p$.  At $t=1$ all coordinates indexed by $\sigma$ vanish, so $H_1(p)\in |C|$.  If $p\in |C|$, then $s(p)=1$ and $H_t(p)=p$ for all $t$.  This proves the strong deformation retraction.

Now suppose that $K$ triangulates $M$ and $\sigma$ is a facet.  The simplex $|\sigma|$ is a PL embedded closed three-ball.  Let
\[
N=M\setminus \operatorname{relint}|\sigma|.
\]
Then $N$ is a compact three-manifold with boundary $\partial N\cong S^2$, and van Kampen applied to
\[
M=N\cup |\sigma|,
\qquad
N\cap|\sigma|=\partial|\sigma|\cong S^2,
\]
gives $\pi_1(N)\cong\pi_1(M)$.  By a collar of $\partial N$, the inclusion $\operatorname{int}N=M\setminus|\sigma|\hookrightarrow N$ is a homotopy equivalence.  The first part now yields
\[
\pi_1(C)\cong\pi_1(M\setminus|\sigma|)\cong\pi_1(M).
\]
\end{proof}

This deformation retraction is also used e.g. in the simplicial-complement arguments of Bagchi and Datta \cite[Section 4]{BagchiDatta}.

\subsection{Coning a disconnected subcomplex}

We shall use a slightly more general form of the Seifert--van Kampen theorem; we include its proof for completeness. 
Let $F_r$ denote the free group of rank $r$.

\begin{lemma}\label{lem:cone}
Let $L$ be a connected CW complex and let
\[
A=A_1\sqcup\cdots\sqcup A_r\subseteq L
\]
be a nonempty subcomplex with $r$ path-connected components.  Attach the cone $x*A$ by its base, and put
\[
X=L\cup_A(x*A).
\]
Choose base points and paths in $L$ so that the images of the groups $\pi_1(A_i)$ are regarded, up to conjugacy, as subgroups of $\pi_1(L)$.  Then
\begin{equation}\label{eq:cone-vk}
\pi_1(X)\cong
\left(
\pi_1(L)
\big/
\normal{\operatorname{im}\pi_1(A_1),\ldots,
\operatorname{im}\pi_1(A_r)}
\right)*F_{r-1}.
\end{equation}
\end{lemma}

\begin{proof}
The groupoid form of the van Kampen theorem applies when the intersection is disconnected \cite[Theorem~3.4]{BrownVanKampen}.
We make the resulting free factor explicit.

For each $i$, take a distinct cone point $x_i$ and form
\[
Y=L\cup_{i=1}^{r}(x_i*A_i).
\]
Since each $A_i$ is connected, successive applications of the ordinary van Kampen theorem give
\[
\pi_1(Y)\cong
\pi_1(L)\big/
\normal{\operatorname{im}\pi_1(A_1),\ldots,
\operatorname{im}\pi_1(A_r)}.
\]
Here the images are identified in $\pi_1(L)$ using the paths chosen in the statement.

For each $i=2,\ldots,r$, attach an edge $e_i$ from $x_1$ to $x_i$, and call the resulting space $Z$.
Because $Y$ is connected, each edge adds one free generator and imposes no relation, so
\[
\pi_1(Z)\cong\pi_1(Y)*F_{r-1}.
\]
The edges $e_2,\ldots,e_r$ and their endpoints form a contractible tree $T\subseteq Z$.
Collapsing $T$ identifies all the cone points and gives $Z/T\cong X$.
Since $T$ is a contractible subcomplex of the CW complex $Z$, the quotient map $Z\to Z/T$ is a homotopy equivalence \cite[Proposition~0.17]{Hatcher}.
Consequently,
\[
\pi_1(X)\cong
\left(
\pi_1(L)\big/
\normal{\operatorname{im}\pi_1(A_1),\ldots,
\operatorname{im}\pi_1(A_r)}
\right)*F_{r-1},
\]
as claimed.
\end{proof}

\subsection{Subcomplexes of a two-sphere}

\begin{lemma}\label{lem:planar}
Let $A$ be a connected finite subcomplex of a triangulated two-sphere.  Then $\pi_1(A)$ is a free group.  In particular,
\[
\pi_1(A)\ne1\quad\Longrightarrow\quad H_1(A;\Q)\ne0,
\]
and if $\dim_{\Q}H_1(A;\Q)=1$, then $\pi_1(A)\cong\Z$.
\end{lemma}

\begin{proof}
A closed regular neighborhood $N$ of $|A|$ in $S^2$ collapses onto $|A|$.  The surface $N$ is connected and has genus zero.  Hence $\pi_1(N)$, and therefore $\pi_1(A)$, is free.  The rank of this free group equals $\dim_{\Q}H_1(A;\Q)$.
\end{proof}

We also use simplicial Alexander duality.  If $S$ is a triangulated two-sphere and $A=S[U]$, $D=S[V(S)\setminus U]$ are complementary induced subcomplexes, then
\begin{equation}\label{eq:alexander}
\widetilde H_1(A;\Q)
\cong
\widetilde H^{0}(D;\Q).
\end{equation}
This is the two-dimensional orientable-field version of \cite[Lemma 4.1]{BagchiDatta}; it also follows directly from classical Alexander duality after replacing the induced subcomplexes by regular neighborhoods.

\subsection{Small complexes and group deficiency}

We isolate the recent computational classification from \cite{GovcEtAl} used in the proof.

\begin{theorem}\label{thm:smallgroups}
If $Y$ is a connected simplicial complex on at most eight vertices, then either $\pi_1(Y)$ is free, or
\[
\pi_1(Y)\cong G*F_s
\]
for some $s\ge0$, where $G$ is one of
\[
\Z/2,
\quad \Z/3,
\quad \Z/4,
\quad \Z^2,
\quad \mathcal K,
\quad B_3.
\]
Here
\[
\mathcal K=\langle a,b\mid aba^{-1}=b^{-1}\rangle
\]
is the Klein-bottle group, and
\[
B_3=\langle a,b\mid aba=bab\rangle
\]
is the braid group on three strands.
\end{theorem}

This is the form of the classification from \cite[Theorem 4.8]{GovcEtAl} needed here, with the free-group case separated from the listed nonfree Grushko factors.

For a finite presentation $\mathcal P=\langle S\mid R\rangle$, write
\[
\defi(\mathcal P)=|S|-|R|.
\]
The deficiency of a finitely presented group $G$ is
\[
\defi(G)=\max_{\mathcal P\text{ presenting }G}\defi(\mathcal P).
\]
For the standard definition and further discussion of group deficiency, see \cite{RobinsonGroups, EpsteinDeficiency, RapaportDeficiency}.

The standard deficiency bound in \cite[Theorem~14.1.5]{RobinsonGroups} also implies that $\defi(\Z^3)\leq 0$.

\begin{lemma}\label{lem:defz3}
The free abelian group $\Z^3$ has deficiency zero.
\end{lemma}

\begin{proof}
The presentation
\[
\Z^3=\langle u,v,w\mid [u,v],[u,w],[v,w]\rangle
\]
has deficiency zero, so $\defi(\Z^3)\ge0$.

Conversely, let
\[
\mathcal P=\langle x_1,\dots,x_g\mid R_1,\dots,R_r\rangle
\]
be any finite presentation of $\Z^3$, and let $X$ be its presentation complex.  Then $X$ is a connected finite two-dimensional CW complex with $\pi_1(X)\cong\Z^3$ and
\[
\chi(X)=1-g+r.
\]
The Hopf exact sequence \cite{Hopf1942} gives a surjection
\[
H_2(X;\Q)\twoheadrightarrow H_2(\pi_1(X);\Q)
=H_2(\Z^3;\Q).
\]
Since $T^3$ is a $K(\Z^3,1)$, the K\"unneth formula gives \cite[Example~3.16]{Hatcher}
\[
H_2(\Z^3;\Q)
\cong H_2(T^3;\Q)
\cong \Lambda^2 H_1(T^3;\Q)
\cong \Lambda^2\Q^3.
\]
Thus $\beta_2(X)\ge\dim_{\Q}\Lambda^2\Q^3=3$.

Since $\beta_1(X)=3$,
\[
1-g+r=\chi(X)
=1-\beta_1(X)+\beta_2(X)
\ge1-3+3=1.
\]
Thus $r\ge g$, so every presentation of $\Z^3$ has deficiency at most zero.  Therefore $\defi(\Z^3)=0$.
\end{proof}

\begin{lemma}\label{lem:small-def}
Let $P$ be the fundamental group of a connected simplicial complex on at most eight vertices.  If $P$ surjects onto $\Z^3$, then $P$ admits a finite presentation of deficiency at least two.
\end{lemma}

\begin{proof}
The surjection $P\twoheadrightarrow\Z^3$ induces a surjection on rational first homology, so
\begin{equation}\label{eq:b1P}
\beta_1(P;\Q)\ge3.
\end{equation}
If $P\cong F_s$ is free, then \eqref{eq:b1P} gives $s\ge3$, and the standard presentation has deficiency $s\ge3$.

Otherwise, Theorem~\ref{thm:smallgroups} gives $P\cong G*F_s$ for one of the listed groups.  Since
\[
\beta_1(P;\Q)=\beta_1(G;\Q)+s,
\]
we obtain the following cases.

\begin{center}
\small
\begin{tabular}{@{}llll@{}}
\toprule
$G$ & $\beta_1(G;\Q)$ & required $s$ & displayed deficiency of $P$ \\
\midrule
$\Z/2,\Z/3,\Z/4$ & $0$ & $s\ge3$ & $s\ge3$ \\
$\Z^2$ & $2$ & $s\ge1$ & $1+s\ge2$ \\
$\mathcal K$ & $1$ & $s\ge2$ & $1+s\ge3$ \\
$B_3$ & $1$ & $s\ge2$ & $1+s\ge3$ \\
\bottomrule
\end{tabular}
\end{center}

Indeed, each finite cyclic group in the table has a one-generator, one-relator presentation; $\Z^2$, $\mathcal K$, and $B_3$ each have a two-generator, one-relator presentation.  Adjoining the free factor $F_s$ adds $s$ generators and no relations.  The minimum lower bound in the table is two.
\end{proof}

\section{The 13-vertex obstruction for $\Tthree$}

\begin{proof}[Proof of Theorem~\ref{thm:main}]
Suppose, for a contradiction, that $K$ is a simplicial triangulation of $\Tthree$ with
\[
f_0(K)=13.
\]

\smallskip
\noindent\textbf{Step 1: the one-skeleton is almost complete.}
Since $\pi_1(K)\cong\Z^3$ requires three generators, \eqref{eq:NS} gives
\[
f_1(K)-42=g_2(K)\ge10\cdot3=30.
\]
Thus
\[
f_1(K)\ge72.
\]
Since a graph on 13 vertices has at most $\binom{13}{2}=78$ edges, the number
\[
q:=78-e
\]
of missing edges satisfies
\begin{equation}\label{eq:q}
q\le6.
\end{equation}

\smallskip
\noindent\textbf{Step 2: a nine-vertex punctured-torus complex.}
Choose a tetrahedral facet
\[
\sigma=\{a,b,c,d\}\in K,
\]
and let 
\[
W=V(K)\setminus\sigma,
\qquad C:=K[W].
\]
By Lemma~\ref{lem:complement},
\begin{equation}\label{eq:piC}
\pi_1(C)\cong\pi_1(\Tthree)\cong\Z^3.
\end{equation}
Moreover, at most the $q\le6$ missing edges of $K$ can be missing from the complete graph on $W$.  Hence
\begin{equation}\label{eq:eC}
f_1(C)\ge\binom92-6=30.
\end{equation}

\smallskip
\noindent\textbf{Step 3: every vertex link in $C$ is connected and non-simply connected.}
Fix $x\in W$ and set
\[
L=C-x,
\qquad
A=\lk_C(x).
\]
Deleting $x$ removes at most eight edges, so \eqref{eq:eC} gives
\[
f_1(L)\ge22.
\]
A disconnected graph on eight vertices has at most $\binom72=21$ edges.  Therefore $L$ is connected.

The star-antistar decomposition at $x$ is
\[
C=L\cup(x*A),
\qquad
L\cap(x*A)=A.
\]
The complex $A$ is nonempty: otherwise $L=C-x$ would contain at least thirty edges on eight vertices, which is impossible.  Suppose that $A$ has $r$ connected components.  By Lemma~\ref{lem:cone},
\[
\pi_1(C)\cong Q*F_{r-1}
\]
for a quotient $Q$ of $\pi_1(L)$.  In view of \eqref{eq:piC}, the left-hand side is $\Z^3$.  If $r\ge2$ and $Q\ne1$, then the right-hand side is a nontrivial free product and is nonabelian.  If $Q=1$, then it is a free group, never isomorphic to $\Z^3$.  Hence $r=1$, thus $A$ is connected.

If $A$ were simply connected, then Lemma~\ref{lem:cone} would give
\[
\pi_1(L)\cong\pi_1(C)\cong\Z^3.
\]
This contradicts Theorem~\ref{thm:smallgroups}, since $L$ has eight vertices and $\Z^3$ is neither free nor one of the listed groups with a free factor.  Therefore
\begin{equation}\label{eq:A-nonsimple}
\lk_C(x)\text{ is connected and non-simply connected for every }x\in W.
\end{equation}

\smallskip
\noindent\textbf{Step 4: a special link has a cyclic fundamental group.}
Consider the triangular face $\tau=\{b,c,d\}$ of $\sigma$.  Since $K$ is a closed combinatorial three-manifold, $\tau$ is contained in exactly two tetrahedra.  One is $\sigma$; write the other as
\[
\{x,b,c,d\}
\]
for a vertex $x\in W$.

Let
\[
S:=\lk_K(x)\cong S^2.
\]
Inside $S$, define the complementary induced subcomplexes
\[
A:=S[W\cap V(S)]=\lk_C(x)
\]
and
\[
D:=S[\sigma\cap V(S)].
\]
By \eqref{eq:A-nonsimple} and Lemma~\ref{lem:planar},
\[
H_1(A;\Q)\ne0.
\]
Simplicial Alexander duality \eqref{eq:alexander} gives
\[
\widetilde H_1(A;\Q)
\cong
\widetilde H^0(D;\Q),
\]
so $D$ is disconnected.

Because $\{x,b,c,d\}$ is a tetrahedron, $D$ contains the full two-simplex with vertex set ${\{b,c,d\}}$.  The only other possible vertex of $D$ is $a$.  Therefore disconnectedness forces 
$D$ to be the disjoint union of a triangle and a point, and Alexander duality then yields
\[
\dim_{\Q}H_1(A;\Q)=1.
\]
By Lemma~\ref{lem:planar},
\begin{equation}\label{eq:piA}
\pi_1(A)\cong\Z.
\end{equation}

\smallskip
\noindent\textbf{Step 5: the deficiency contradiction.}
Return to
\[
C=L\cup(x*A),
\qquad L=C-x,
\]
for the special vertex $x$ selected in Step 4.  Since $A$ is connected, Lemma~\ref{lem:cone} and \eqref{eq:piA} imply that, for some $w\in \pi_1(L)$,
\begin{equation}\label{eq:one-relator}
\Z^3\cong\pi_1(C) \cong  \pi_1(L)/\normal{\text{im}(\pi_1(A))}
=
\pi_1(L)/\normal{w}.
\end{equation}
In particular, $\pi_1(L)$ surjects onto $\Z^3$.  The complex $L$ has eight vertices, so Lemma~\ref{lem:small-def} gives a presentation
\[
\pi_1(L)=\langle y_1,\dots,y_g\mid R_1,\dots,R_r\rangle
\]
with
\[
g-r\ge2.
\]
Choose a word $W=W(y_1,\dots,y_g)$ representing $w$.  Equation \eqref{eq:one-relator} gives a presentation
\[
\Z^3\cong
\langle y_1,\dots,y_g\mid R_1,\dots,R_r,W\rangle
\]
of deficiency
\[
g-(r+1)\ge1.
\]
This contradicts Lemma~\ref{lem:defz3}, which states that $\defi(\Z^3)=0$.  The contradiction completes the proof.
\end{proof}

\begin{proof}[Proof of Corollary~\ref{cor:bounds}]
Let $K$ triangulate $\Tthree$ with $n=f_0(K)$.  By \eqref{eq:NS},
\[
30\le g_2(K)
=f_1(K)-4n+10
\le\binom n2-4n+10
=\binom{n-4}{2}.
\]
For $n\le12$, the right-hand side is at most $\binom82=28$, a contradiction.  Theorem~\ref{thm:main} excludes $n=13$, and the K\"uhnel--Lassmann construction gives $n\le15$.
\end{proof}

\end{document}